\documentclass[12pt, a4paper]{article}
\usepackage{amscd}
\newtheorem{theorem}{Theorem}

\newtheorem{proposition}[theorem]{Proposition}
\newtheorem{remark}[theorem]{Remark}

\def\R{I\kern-.3em R}
\def\P{I\kern-.3em P}
\def\N{I\kern-.3em N}
\def\Qed{\hfill\raisebox{.6ex}{\framebox[2.5mm]{}}\\[.15in]}
\title{On surfaces with $p_g=q=1$ and non-ruled bicanonical involution}
\author{Carlos Rito}
\date{}
\begin{document}
\maketitle

\begin{abstract}
This paper classifies surfaces of general type $S$ with $p_g=q=1$
having an involution $i$ such that $S/i$ has non-negative Kodaira
dimension and that the bicanonical map of $S$ factors through
the double cover induced by $i.$

It is shown that $S/i$ is regular and either: a) the Albanese
fibration of $S$ is of genus 2 or b) $S$ has no genus 2 fibration
and $S/i$ is birational to a $K3$ surface. For case a) a list of
possibilities and examples are given. An example for case b) with
$K^2=6$ is also constructed.

\noindent 2000 Mathematics Classification: 14J29.
\end{abstract}
\section{Introduction}
Let $S$ be a smooth irreducible projective surface of general type.
The {\em pluricanonical map} $\phi_n$ of $S$ is the map given by the
linear system $|nK_S|,$ where $K_S$ is the canonical divisor of $S.$
For minimal surfaces $S,$ $\phi_n$ is a birational morphism if
$n\geq 5$ (cf. \cite [Ch.VII, Theorem  (5.2)] {BPV}).
The {\em bicanonical map}
$$\phi_2:S\longrightarrow\P^{K_S^2+\chi(S)-1}$$ is a morphism if
$p_g(S)\geq 1$ (this result is due to various authors, see \cite{CC}
for more details). This paper focuses on the study of surfaces
$S$ with $p_g(S)=q(S)=1$ having an involution $i$ such that the
Kodaira dimension of $S/i$ is non-negative and $\phi_2$ is composed
with $i,$ i.~e. it factors through the double cover $p:S\rightarrow
S/i.$

There is an instance where the bicanonical map is necessarily
composed with an involution: suppose that $S$ has a fibration of
genus 2, i.~e. it has a morphism $f$ from $S$ to a curve such that a
general fibre $F$ of $f$ is irreducible of genus 2. The system
$|2K_S|$ cuts out on $F$ a subseries of the bicanonical series of
$F,$ which is composed with the hyperelliptic involution of $F,$ and
then $\phi_2$ is composed with an involution. This is the so called
{\sl standard case} of non-birationality of the bicanonical map.

By the results of Bombieri, \cite{Bomb}, improved later by
Reider, \cite{Rd}, a  minimal surface $S$ satisfying  $K^2>9$ and
$\phi_2$ non-birational necessarily presents the  standard case of
non-birationality of the bicanonical map.

The non-standard case of non-birationality of the bicanonical  map,
i.~e. the case where $\phi_2$ is non-birational and the surface has
no genus 2 fibration, has been studied by several authors.

Du Val, \cite{Du}, classified the regular surfaces $S$ of general type
with $p_g\geq 3$, whose general canonical curve is smooth
and hyperelliptic.  Of course that, for these surfaces, the bicanonical map is
composed
with an involution $i$ such that $S/i$ is rational.
The families of surfaces exhibited by Du Val, presenting the non-standard case,
are nowadays called the {\em Du Val examples}.

Other authors have later studied the non-standard case: the
articles \cite{CCM}, \cite{CFM}, \cite{CM}, \cite{CM2}, \cite{X} and
\cite{Bo} treat the cases $\chi(\mathcal{O}_S)>1$ or $q(S)\geq 2$
(cf. the expository paper \cite{Ci} for more information on this
problem).

Xiao Gang,  \cite{X}, presented a list of possibilities for the
non-standard case of non-birationality of the  bicanonical morphism
$\phi_2$. For the case when  $\phi_2$ has  degree 2 and the bicanonical
image is a ruled surface, Theorem 2 of \cite{X} extended Du Val's list to
$p_g(S)\geq 1$ and added two extra families (this result is still valid
assuming only that $\phi_2$ is composed with an involution such that
the quotient surface is a ruled surface).
Recently G. Borrelli \cite{Bo} excluded these two families,
confirming that the only possibilities for this instance are the Du Val examples.

For irregular surfaces  the following holds (see
\cite [Theorems 1, 3]{X}, \cite [Theorem A] {CCM},
 \cite [Theorem 1.1] {CM}, \cite {CM2}):\\\\
{\em Suppose that $S$ is a smooth minimal irregular surface of
general type having non-birational bicanonical map.  If $p_g(S)\geq
2$ and $S$ has no genus 2 fibration, then only the following
(effective) possibilities occur:
\begin{itemize}
     \item [$\cdot$] $p_g(S)=q(S)=2,$ $K_S^2=4;$
     \item [$\cdot$] $p_g(S)=q(S)=3,$ $K_S^2=6.$
\end{itemize}
In both cases $\phi_2$ is composed with an involution $i$ such that ${\rm Kod}(S/i)=2$. }\\\\

This paper completes this result classifying the minimal surfaces
$S$  with $p_g(S)=q(S)=1$ such that  $\phi_2$ is composed with an
involution $i$ satisfying ${\rm Kod}(S/i)\geq 0$.

The main result is the following:
\begin{theorem}\label{next}
Let $S$ be a smooth minimal irregular surface of general type with
an involution $i$ such that ${\rm Kod}(S/i)\geq 0$ and the
bicanonical map $\phi_2$ of $S$ is composed with $i.$ If
$p_g(S)=q(S)=1,$ then only the following possibilities can occur:
\begin{itemize}
  \item [{\rm a)}] $S/i$ is regular, the Albanese fibration of $S$ has genus 2 and
\begin{itemize}
   \item [{\rm (i)}] ${\rm Kod}(S/i)=2,$ $\chi(S/i)=2,$ $K_S^2=2,$ ${\rm deg}
   (\phi_2)=8,$ or
   \item [{\rm (ii)}] ${\rm Kod}(S/i)=1,$ $\chi(S/i)=2,$ $2\leq K_S^2\leq 4,$
   ${\rm deg}(\phi_2)\geq 4,$ or
   \item [{\rm (iii)}] $S/i$ is birational to a $K3$ surface, $3\leq K_S^2\leq
   6,$ ${\rm deg}(\phi_2)=4;$
\end{itemize}
  \item [{\rm b)}] $S$ has no genus 2 fibration and $S/i$ is birational to a $K3$ surface.
\end{itemize}

Moreover, there are examples for {\rm(i)}, {\rm (ii)} with $K_S^2=4,$
{\rm (iii)} with $K_S^2=3,4$ or 5  and for {\rm b)} with $K_S^2=6$
and $\phi_2$ of degree 2.
\end{theorem}
\begin{remark}
Examples for {\rm (iii)} were given by Catanese in \cite{C2}. The
other examples will be presented in Section \ref{Examples}.
\end{remark}

Note that surfaces of general type with $p_g=q=1$ and $K^2=3$ or 8 were also
studied by Polizzi in \cite{P} and \cite{P2}.

In the example in Section \ref{Examples} for case b) of Theorem
\ref{next}, $S$ has $p_g=q=1$ and $K^2=6.$  This seems to be the
first construction of a surface with these invariants. This example
contradicts a result of Xiao Gang. More precisely, the list of
possibilities in \cite{X} rules out the case where $S$ has no genus
2 fibration, $p_g(S)=q(S)=1$ and $S/i$ is birational to a $K3$
surface. In Lemma 7 of \cite{X} it is written that $R$ has only
negligible singularities, but the possibility
$\chi(K_{\widetilde{P}}+\widetilde{\delta})<0$ in formula (3) of
page 727 was overlooked. In fact we will see that $R$
($\overline{B}$ in our notation) can have a non-negligible
singularity.

An important technical tool that will be used several times is the
{\em canonical resolution} of singularities of a surface. This is a
resolution of singularities as described in \cite{BPV}.

The paper is organized as follows. Section \ref{Involution} studies
some general properties of a surface of general type $S$ with an
involution $i$. Section \ref{q1} states some properties of surfaces
with $p_g=q=1$. Section \ref{Classification} contains the proof of
Theorem \ref{next}. Crucial ingredients for this proof are the
existence of the Albanese fibration of $S$ and the formulas of
Section \ref{Involution}. In Section \ref{Examples} examples for
Theorem \ref{next} are obtained, via the construction of branch
curves with appropriate singularities. The Computational Algebra
System {\em Magma} is used to perform the necessary calculations
(visit http://magma.maths.usyd.edu.au/magma for more information
about Magma).

\medskip
\noindent{\bf Acknowledgements.} First of all I thank my supervisor
Margarida Mendes Lopes for all her support. I would also like to thank
Miles Reid for the useful conversations and hospitality during my one
month stay in Warwick, in July 2004, where I started learning the
Computational System Magma and where part of this work began.
This stay was partially supported by the european contract EAGER no.
HPRN-CT-2000-00099 and the University of Warwick. I also thank the
Departamento de Matem\' atica da Universidade de Tr\'as-os-Montes e
Alto Douro.

The author is a collaborator of the Center for Mathematical
Analysis, Geometry and Dynamical Systems of Instituto Superior T\'
ecnico, Universidade T\' ecnica de Lisboa and was partially
supported by FCT (Portugal) through program POCTI/FEDER and Project
POCTI/MAT/44068/2002.

\medskip
\noindent{\bf Notation and conventions.} We work over the complex
numbers; all varieties are assumed to be projective algebraic.  We
do not distinguish between line bundles and divisors on a smooth
variety. Linear equivalence is denoted by $\equiv$. A {\em nodal
curve} or {\em $(-2)$-curve} $C$ on a surface is a curve isomorphic
to $\P^1$ such that $C^2=-2$. Given a surface $X,$ ${\rm Kod}(X)$
means the {\em Kodaira dimension} of $X.$ We say that a curve
singularity is {\em negligible} if it is either a double point or a
triple point which resolves to at most a double point after one
blow-up. A $(n,n)$ {\em point}, or {\em point of type} $(n,n),$ is a
point of multiplicity $n$ with an infinitely near point also of
multiplicity $n.$ An {\em involution} of a surface $S$ is an
automorphism of $S$ of order 2. We say that a map is composed with
an involution $i$ of $S$ if it factors through the map $S\rightarrow
S/i.$ The rest of the notation is standard in Algebraic Geometry.
\section{Generalities on involutions}\label{Involution}
Let $S$ be a smooth minimal surface of general type with an
involution $i.$ As $S$ is minimal of general type, this
involution is biregular. The fixed locus of $i$ is the union of a
smooth curve $R''$ (possibly empty) and of $t\geq 0$ isolated points
$P_1,\ldots,P_t.$ Let $S/i$ be the quotient of $S$ by $i$ and
$p:S\rightarrow S/i$ be the projection onto the quotient. The
surface $S/i$ has nodes at the points $Q_i:=p(P_i),$ $i=1,\ldots,t,$
and is smooth elsewhere. If $R''\not=\emptyset,$ the image via $p$
of $R''$ is a smooth curve $B''$ not containing the singular points
$Q_i,$ $i=1,\ldots,t.$ Let now $h:V\rightarrow S$ be the blow-up of
$S$ at $P_1,\ldots,P_t$ and set $R'=h^*R''$. The involution $i$
induces a biregular involution $\widetilde{i}$ on $V$ whose fixed
locus is $R:=R'+\sum_1^t h^{-1}(P_i).$ The quotient
$W:=V/\widetilde{i}$ is smooth and one has a commutative diagram:
$$
\begin{CD}\ V@>h>>S\\ @V\pi VV  @VV p V\\ W@>g >> S/i
\end{CD}
$$
where $\pi:V\rightarrow W$ is the projection onto the quotient and
$g:W\rightarrow S/i$ is the minimal desingularization map. Notice
that $$A_i:=g^{-1}(Q_i),\ \ i=1,\ldots,t,$$ are $(-2)$-curves and
$\pi^*(A_i)=2\cdot h^{-1}(P_i).$ Set $B':=g^*(B'').$ Because $\pi$ is
a double cover with branch locus $B'+\sum_1^t A_i,$ there exists
a line bundle $L$ on $W$ such that $$2L\equiv B:=B'+\sum_1^t A_i.$$
It is well known that (cf. \cite[Ch. V, \S 22]{BPV}):\\
$$
\begin{array}{c}
  p_g(S)=p_g(V)=p_g(W)+h^0(W,\mathcal{O}_W(K_W+L)), \\\\
  q(S)=q(V)=q(W)+h^1(W,\mathcal{O}_W(K_W+L))
\end{array}
$$
and
\begin{equation}\label{relate}
\begin{array}{c}
  K_S^2-t=K_V^2=2(K_W+L)^2, \\\\
  \chi(\mathcal{O}_S)=\chi(\mathcal{O}_V)=2\chi(\mathcal{O}_W)+\frac{1}{2}L(K_W+L).
\end{array}
\end{equation}
Furthermore, from the papers \cite{CM} and \cite{CCM1}, if $S$ is a smooth minimal
surface of general type with an involution $i,$ then
\begin{equation}\label{Mrg1}
\chi(\mathcal{O}_W(2K_W+L))=h^0(W,\mathcal{O}_W(2K_W+L)),
\end{equation}
\begin{equation}\label{Mrg2}
\chi(\mathcal{O}_W)-\chi(\mathcal{O}_S)=K_W(K_W+L)-h^0(W,\mathcal{O}_W(2K_W+L))
\end{equation}
and the bicanonical map
\begin{equation}\label{Mrg3}
 \phi_2 {\rm\ is\ composed\ with\ } i {\rm\ if\ and\ only\ if\ } h^0(W,\mathcal{O}_W(2K_W+L))=0.
\end{equation}
From formulas (\ref{relate}) and (\ref{Mrg2}) one obtains
the number $t$ of nodes of $S/i:$
\begin{equation}\label{numbernodes}
t=K_S^2+6\chi(\mathcal{O}_W)-2\chi(\mathcal{O}_S)-2h^0(W,\mathcal{O}_W(2K_W+L)).
\end{equation}

Let $P$ be a minimal model of the resolution $W$ of $S/i$ and
$\rho:W\rightarrow P$ be the natural projection. Denote by
$\overline{B}$ the projection $\rho(B)$ and by $\delta$ the "projection" of $L.$
\begin{remark}\label{CanRes}
Resolving the singularities of $\overline{B}$ we obtain exceptional
divisors $E_i$ and numbers $r_i\in 2\N^+$ such that $E_i^2=-1,$
$K_W=\rho^*(K_P)+\sum E_i$ and $B=\rho^*(\overline{B})-\sum r_iE_i.$
\end{remark}
\begin{proposition}\label{equation}
With the previous notations, the bicanonical map $\phi_2$ is composed with $i$ if and only if
$$\chi(\mathcal{O}_P)-\chi(\mathcal{O}_S)=
  K_P(K_P+\delta)+\frac{1}{2}\sum(r_i-2).$$
\end{proposition}
{\bf Proof:}
  From formulas (\ref{Mrg2}), (\ref{Mrg3}) and Remark \ref{CanRes} we get
  $$\chi(\mathcal{O}_P)-\chi(\mathcal{O}_S)=\frac{1}{2}K_W(2K_W+2L)=$$
  $$=\frac{1}{2}\left(\rho^*(K_P)+\sum E_i\right)
  \left(2\rho^*(K_P+\delta)+\sum(2-r_i)E_i\right)=$$
  $$=K_P(K_P+\delta)+\frac{1}{2}\sum(r_i-2).$$ \Qed
\section{Surfaces with $p_g=q=1$ and an involution}\label{q1}

Let $S$ be a minimal smooth projective surface of general type satisfying $p_g(S)=q(S)=1$.

Note that then $2\leq K_S^2\leq 9$: we have $K_S^2\leq
9\chi(\mathcal{O}_S)$ by the  Myiaoka- -Yau inequality (see
\cite[Ch. VII, Theorem (4.1)]{BPV})  and $K_S^2\geq 2p_g$ for an
irregular surface (see \cite{De}).

Furthermore, if the bicanonical map of $S$ is not birational, then $K_S^2\neq 9$.
In fact, by \cite{CM}, if $K_S^2=9$ and
$\phi_2$ is not birational, then $S$ has a genus 2 fibration, while
Th\' eor\` eme 2.2 of \cite{X1} implies that if $S$ has a genus 2
fibration and $p_g(S)=q(S)=1,$ then $K_S^2\leq 6.$

Since $q(S)=1$ the Albanese variety of $S$ is an elliptic curve $E$ and the Albanese
map is a connected fibration (see e.~g. \cite{B} or \cite{BPV}).

Suppose that $S$ has an involution $i.$ Then $i$ preserves the Albanese fibration
(because $q(S)=1$) and so we have a commutative diagram
\begin{equation}\label{CmtDgr}
\begin{CD}\ V@>h>>S @>>>E\\ @V\pi VV  @VV p V @VV V\\ W@>>> S/i @>>>\Delta
\end{CD}
\end{equation}
where $\Delta$ is a curve of genus $\leq 1.$ Denote by
$$f_A:W\rightarrow\Delta$$ the fibration induced by the Albanese
fibration of $S.$

Recall that $$\rho:W\rightarrow P$$ is the projection of $W$ onto
its minimal model $P$ and $$\overline{B}:=\rho(B),$$ where
$B:=B'+\sum_1^t A_i\subset W$ is the branch locus of $\pi.$

Let $$\overline{B'}:=\rho(B'),\ \ \overline{A_i}=\rho(A_i).$$

When $\overline{B}$ has only negligible singularities, the map $\rho$
contracts only exceptional curves contained in fibres of $f_A.$ In
fact, there exists otherwise a $(-1)$-curve $J\subset W$ such that
$JB=2$ and so $\pi^*(J)$ is a rational curve transverse to the
fibres of the (genus 1 base) Albanese fibration of $S,$ which is
impossible. 
Moreover, $\rho$ contracts no curve meeting $\sum A_i,$
because $h:V\rightarrow S$ is the contraction of isolated $(-1)$-curves.
Therefore the singularities of $\overline{B}$ are exactly the
singularities of $\overline{B'},$ i.~e. $\overline{B'}\bigcap\sum
\overline{A_i}=\emptyset.$ In this case the image of $f_A$ on $P$
will be denoted by $\overline{f_A}.$

If $\Delta\cong\P^1,$ then the double cover $E\rightarrow\Delta$ is
ramified over 4 points $p_j$ of $\Delta,$ thus the branch locus
$B'+\sum_1^t A_i$ is contained in 4 fibres $$F_A^j:=f_A^*(p_j),\ j=1,...,4,$$ of the fibration $f_A$.
Hence, by Zariski's Lemma (see e. g. \cite{BPV}), the irreducible
components $B_i'$ of $B'$ satisfy $B_i'^2\leq 0.$ If $\overline{B}$
has only negligible singularities, then also $\overline{B'}^2\leq 0.$
As  $\pi^*(F_A^j)$ is of multiplicity 2, each component of
$F_A^j$ which is not a component of the branch locus $B'+\sum_1^t
A_i$ must be of even multiplicity.
\section{The classification theorem}\label{Classification}
In this section we will prove Theorem \ref{next}. We will use freely the notation and results
of Sections \ref{Involution} and \ref{q1}.\\\\
{\bf Proof of Theorem \ref{next}:} Since $p_g(P)\leq p_g(S)=1,$ then
$\chi(\mathcal{O}_P)\leq 2-q(P)\leq 2.$ Proposition \ref{equation}
gives $\chi(\mathcal{O}_P)\geq 1,$ because $K_P$ is nef (i.~e. $K_PC\geq 0$
for every curve $C$). So from Proposition \ref{equation} and the classification
of surfaces (see e.~g. \cite{B} or \cite{BPV}) only the following cases can occur:
\begin{itemize}
     \item [1)] $P$ is of general type;
     \item [2)] $P$ is a surface with Kodaira dimension 1;
     \item [3)] $P$ is an Enriques surface, $\overline{B}$ has only
   negligible singularities;
     \item [4)] $P$ is a K3 surface, $\overline{B}$ has a 4-uple or
   $(3,3)$ point, and possibly negligible singularities.
\end{itemize}
We will show that: case 3) does not occur, in cases 1) and 2) the
Albanese fibration has genus 2 and  only in case 4)  the Albanese
fibration can have genus $\neq 2$.

Each of cases 1),..., 4) will be studied separately. We start by considering:
\medskip

     {\bf Case 1)} As $P$ is of general type, $K_P^2\geq 1$ and $K_P$ is nef,
     Proposition \ref{equation} gives $\chi(\mathcal{O}_P)=2,$
     $K_P^2=1,$ $K_P\delta=0$ and $\overline{B}$ has only negligible
     singularities. The equality $K_P\overline{B'}=K_P2\delta=0$
     implies $\overline{B'}^2<0$ when $B'\ne 0.$
     In the notation of Remark \ref{CanRes} one has
     $K_W\equiv \rho^*(K_P)+\sum E_i$ and $B'=\rho^*(\overline{B'})-2\sum E_i.$
     So
     $$K_S^2=K_V^2+t=\frac{1}{4}(2K_V)^2+t=\frac{1}{4}\pi^*(2K_W+B)^2+t=$$
     $$=\frac{1}{2} (2K_W+B)^2+t=\frac{1}{2}(2K_W+B')^2=\frac{1}{2}(2K_P+\overline{B'})^2=
     \frac{1}{2}(4+\overline{B'}^2).$$
     Since $K_S^2\geq 2p_g(S)$ for an irregular surface (see \cite{De}),
     $\overline{B'}^2<0$ is impossible, hence $B'=0$ and $K_S^2=2.$
     By \cite{C} minimal surfaces of general type with $p_g=q=1$
     and $K^2=2$ have Albanese fibration of genus 2.
     This is case (i) of Theorem \ref{next}.
     We will see in Section \ref{Examples} an example for this case.

     Finally the fact that ${\rm deg}(\phi_2)=8$ follows immediately because
     $\phi_2$ is a morphism onto $\P^2$ and $(2K_S)^2=8$.\\

Next we exclude:
\medskip

     {\bf Case 3)} Using the notation of Remark \ref{CanRes} of Section \ref{Involution}, we can write
     $K_W\equiv\rho^*(K_P)+\sum E_i$ and $2L\equiv \rho^*(2\delta)-2\sum E_i,$
     for some exceptional divisors $E_i.$
     Hence $$L(K_W+L)=\frac{1}{2}L(2K_W+2L)=$$
     $$=\frac{1}{2}(\rho^*(\delta)-\sum E_i)(2\rho^*(K_P)+\rho^*(2\delta))=
     \frac{1}{2}\delta(2K_P+2\delta)=\delta^2$$
     and then, from (\ref{relate}), $\delta^2=-2.$
     Now (\ref{Mrg3}) and (\ref{numbernodes}) imply $t=K_S^2+4,$
     thus $$\overline{B'}^2=\overline{B}^2+2t=(2\delta)^2+2t=-8+2t=2K_S^2>0.$$
     This is a contradiction because we have seen that $\overline{B'}^2\leq 0$
     when $\overline{B'}$ has only negligible singularities.
     Thus case 3) does not occur.\\

Now we focus on:
\medskip

{\bf Case 2)}
     Since we are assuming that  ${\rm Kod}(P)=1$, $P$ has an elliptic fibration
     (i.~e. a morphism $f_e: P\to C$ where $C$ is a curve and the general fibre of $f_e$
     is a smooth connected elliptic curve). Then $K_P$ is
     numerically equivalent to a rational multiple of a fibre of $f_e$
     (see e.~g. \cite{B} or \cite{BPV}). As $K_P\delta\geq 0$,
      Proposition \ref{equation}, together with
     $\chi(\mathcal{O}_P)\leq 2,$  yield $K_P\delta=0$ or 1.

     Denote by $F_e$ (resp. $F_A$) a general fibre of $f_e$ (resp. $f_A$) and let
     $\overline{F_A}:=\rho(F_A).$
     If $K_P\delta=0,$ then $F_e\overline{B}=0,$ which implies that
     the fibration $f_e$ lifts to an elliptic fibration on  $S.$ This is
     impossible because $S$ is a surface of general type. So $K_P\delta=1$
     and, since $p_g(P)\leq p_g(S)=1,$ the only possibility allowed by
     Proposition \ref{equation} is
\begin{itemize}
   \item [] $p_g(P)=1,$ $q(P)=0$ and
     $\overline{B}$ has only negligible singularities.
\end{itemize}
     Now $q(P)=0$ implies that the elliptic fibration $f_e$ has a rational base,
     thus the canonical bundle formula (see e.~g. \cite[Ch. V, \S 12]{BPV}) gives
     $K_P\equiv\sum(m_i-1)F_i,$ where $m_iF_i$ are the multiple fibres of $f_e.$
     From
     $$2=2\delta K_P=\overline{B'}K_P=\overline{B'}\sum(m_i-1)F_i,
     \ \ \ \ \ \overline{B'}F_i\equiv 0\ ({\rm mod\ 2})$$ we get
     $$K_P\equiv\frac{1}{2}F_e.$$

     Since $\overline{B}$ has only negligible singularities, $\overline{B'}^2\leq 0$
     and then
     \begin{equation}\label{eqKB}
     2K_S^2=(2K_W+B')^2=\rho^*\left(2K_P+\overline{B'}\right)^2=8+\overline{B'}^2\leq 8.
     \end{equation}

     Therefore $2\leq K_S^2\leq 4$.
     If $K_S^2=2,$ then the Albanese fibration of $S$
     is of genus 2, by \cite{C}.
     So, to prove  statement  a),\ (ii) of Theorem \ref{next}, we must show that for
     $K_S^2=3$ or 4 the Albanese fibration of $S$ has genus 2.
     We will study each of these cases separately.\\\\
     First we consider\\\\
     $\cdot$ ${\bf K_S^2=4.}$\\\\
     Let $\overline{F_A^i}:=\rho(F_A^i),$ $i=1,\ldots,4.$
\\\\
{\bf Claim 1:} {\em If $f_A$ is not a genus 2 fibration then
     $$\overline{F_A^j}=2\overline{B'},$$} for some $j\in\{1,\ldots,4\}.$
     \\\\
     {\em Proof}\ : By formula (\ref{eqKB})  $\overline{B'}^2=0,$ and so $\overline{B'}$
     contains the support of $x\geq 1$ of the $\overline{F_A^i}$'s. The facts
     $K_P\overline{F_A}>0$
     (because $g(\overline{F_A})\geq 2$) and $K_P\overline{B'}=2$ imply $x=1,$ i.~e.
     $\overline{F_A^j}=k\overline{B'},$ for some $j\in\{1,\ldots,4\}$ and $k\in\N^+.$
     If $k=1$ then $\overline{F_A}K_P=2,$ thus $\overline{F_A}$ is of genus 2
     and $S$ is as in case (ii) of Theorem \ref{next}.

     Suppose now $k\geq 2.$ Then each irreducible component of the divisor
     $$D:=\overline{F_A^1}+\ldots+\overline{F_A^4}$$ whose support is not
     in $\sum_1^{14}\overline{A_i}$ is of multiplicity greater than 1.
     The fibration $\overline{f_A}$ gives a cover $F_e\rightarrow\P^1$ of degree
     $\overline{F_A}F_e,$ for a general fibre $F_e$ of the elliptic fibration $f_e$.
     The Hurwitz formula (see e.~g. \cite{GH}) says that the ramification degree $r$
     of this cover is $2\overline{F_A}F_e.$
     Let $p_1,\ldots,p_n$ be the points in $F_e\cap D$
     and $\alpha_i$ be the intersection number of $F_e$ and $D$ at
     $p_i.$ Of course $F_eD=4\overline{F_A}F_e=\sum_1^n \alpha_i$ and then
     $F_e\bigcap\sum \overline{A_i}=\emptyset$ implies $\alpha_i\geq
     2,$ $i=1,\ldots,n.$
     We have $$2\overline{F_A}F_e=r\geq\sum_1^n(\alpha_i-1)=\sum_1^n\alpha_i-n=4\overline{F_A}F_e-n,$$
     i.~e. $n\geq 2\overline{F_A}F_e.$
     The only possibility is $n=2\overline{F_A}F_e$ and $\alpha_i=2\ \forall i,$ which means
     that every component $\Gamma$ of $D$ such that $\Gamma F_e\not=0$ is
     exactly of multiplicity 2. In particular an irreducible component of $\overline{B'}$
     is of multiplicity 2, thus $k=2,$ i.~e. $\overline{F_A^j}=2\overline{B'}.$ $\diamond$
     \\\\
     {\bf Claim 2:} {\em There is a smooth rational curve $C$ contained in a
     fibre $F_C$ of the elliptic fibration $f_e,$ and not contained in fibres of
     $\overline{f_A},$ such that
     \begin{equation}\label{eqC3}
     m:=\widehat{C}\sum_1^t A_i\leq 3,
     \end{equation}
     where $\widehat{C}$ is the strict transform of
     $C$ in $W.$}
     \\\\
     {\em Proof}\ : Since $\overline{A_i}F_e=\overline{A_i}2K_P=0,$
     then each $\overline{A_i}$
     is contained in a fibre of $f_e,$ and in particular the elliptic fibration
     $f_e$ has reducible fibres.
     Denote by $C$ an irreducible component of a reducible fibre $F_C$ of $f_e,$
     by $\xi$ the multiplicity of $C$ in $F_C$
     and by $\widehat{C}$ the strict transform of $C$ in $W.$
     If the intersection number of $C$ and the support of $F_C-\xi C$ is
     greater than 3 then, from the configurations of
     singular fibres of an elliptic fibration (see e.~g. \cite[Ch. V, \S 7]{BPV}),
     $F_C$ must be   of type $I_0^*,$ i.~e. it has the following configuration:
     it is the union of four disjoint $(-2)$-curves $\theta_1,\ldots,\theta_4$
     with a $(-2)$-curve $\theta,$ with multiplicity 2, such that
     $\theta \theta_i=1,$ $i=1,\ldots,4.$

     So if $\widehat C\sum_1^tA_i>3,$ the fibre $F_C$ containing $C$ is of type $I_0^*$
     with $\widehat C\sum_1^tA_i=4.$ Since the number of nodes of $S/i$ is 
     $t=K_S^2+10=14\not\equiv 0\ {\rm (mod\ 4)},$ there must be a reducible fibre such that 
     for every component $C\not\subset\sum_1^t\overline{A_i},$ $\widehat C\sum_1^tA_i\leq 3.$
     As $f_e\ne\overline{f_A}$ and the $\overline{A_i}$'s are
     contained in fibres of $f_e$ and in fibres of $\overline{f_A},$ we can choose $C$
     not contained in fibres of $\overline{f_A}.$\ $\diamond$\\

     Let $C$ be as in Claim 2 and consider the resolution
     $\widetilde{V}\rightarrow V$ of the
     singularities of $\pi^*(\widehat{C}).$ Let $G\subset \widetilde{V}$
     be the strict transform of $\pi^*(\widehat{C}).$
     Notice that $G$ has multiplicity 1, because $C$ transverse to the fibres of
     $\overline{f_A}$ implies $C\not\subset \overline{B}.$ Recall that $E$
     denotes the basis of the Albanese fibration of $S.$
     \\\\
     {\bf Claim 3:} {\em The Albanese fibration of $\widetilde{V}$ induces a
     cover $G\rightarrow E$ with ramification degree $$r:=K_{\widetilde{V}}G+G^2.$$}
     \\\\
     {\em Proof}\ : Let $G_1,\ldots,G_h$ be the connected (hence smooth)
     components of $G.$
     The curve $C$ is not contained in fibres of $\overline{f_A},$ thus $G$ is
     not contained in fibres of the Albanese fibration of $\widetilde{V}$.
     This fibration induces a cover $G_i\rightarrow E$ with ramification degree,
     from the Hurwitz formula,
     $$r_i=2g(G_i)-2=K_{\widetilde{V}}G_i+G_i^2.$$
     This way we have a cover $G\rightarrow E$ with ramification degree
     $$r=\sum r_i=K_{\widetilde{V}}(G_1+\cdots+G_h)+\left(G_1^2+\cdots+G_h^2\right)=
     K_{\widetilde{V}}G+G^2.\ \diamond$$\medskip

     We are finally in position to show that $g(F_A)=2.$

     Let $n:=\widehat{C}B'.$ We have
     $$2K_V\pi^*(\widehat{C})=\pi^*(2K_W+B'+\sum A_i)\pi^*(\widehat{C})=$$
     $$=2(2K_W+B'+\sum A_i)\widehat{C}=4K_W\widehat{C}+2(B'+\sum A_i)\widehat{C}=$$
     $$=4(-2-\widehat{C}^2)+2(n+m)=-8-2\pi^*(\widehat{C})^2+2(n+m),$$
     i.~e. $$K_V\pi^*(\widehat{C})+\pi^*(\widehat{C})^2=n+m-4.$$

     Suppose that $g(F_A)\ne 2.$
     Let $\Lambda\subset V$ be the double Albanese fibre induced by
     $\overline{F_A^j}=2\overline{B'}$ (as in Claim 1) and $\widetilde{\Lambda}\subset \widetilde{V}$
     be the total transform of $\Lambda.$ From
     $$G\widetilde{\Lambda}=\pi^*(\widehat{C})\Lambda\geq\pi^*(\widehat{C})\pi^*(B')=2n$$
     one has $r\geq n.$
     Then
     $$n+m-4=K_V\pi^*(\widehat{C})+\pi^*(\widehat{C})^2\geq K_{\widetilde{V}}G+G^2=r\geq n$$
     and so $m\geq 4,$ which contradicts Claim 2.

     So if $K_S^2=4,$ then the Albanese fibration of $S$ is of genus 2.
     \\\\
     We will now consider the possibility
     \\\\
     $\cdot$ ${\bf K_S^2=3.}$
     \\\\
     In this case a general Albanese
     fibre $\Lambda$ has genus 2 or 3 (see \cite{CC}). Suppose then $g(\Lambda)=3.$
     Surfaces $S$ with $K_S^2=g(\Lambda)=3$ are studied in detail in \cite{CC}.
     There (see also \cite{K}) it is shown that the relative canonical map
     $\gamma,$ given by $|K_S+n\Lambda|$ for some $n,$ is a morphism.

     We know that $K_P\overline{B'}=2$ and $\overline{B'}^2=-2,$ by
     (\ref{eqKB}). We have already seen that $\overline{B}$ has only
     negligible singularities (which means $r_i=2\ \forall i,$ in the
     notation of Remark \ref{CanRes}) and then $\rho$ contracts no curve
     meeting $\sum A_i.$
     Let $R'$ be the support of $\pi^*(B').$
     \\\\
     {\bf Claim 4:} {\em We have $$K_VR'=1.$$}

     \noindent {\em Proof}\ :
     $$2K_V\cdot 2R'=\pi^*(2K_W+B)\pi^*(B')=2(2K_W+B)B'=$$
     $$=2(2K_W+B')B'=2\left(2\rho^*(K_P)+\rho^*(\overline{B'})\right)
     \left(\rho^*(\overline{B'})-\sum 2E_i\right)=$$
     $$=2(2K_P+\overline{B'})\overline{B'}=2(4-2)=4,$$
     thus $K_VR'=1.$ $\diamond$
     \\

     As the map $$\gamma\circ h:V\longrightarrow\gamma(S)$$
     is a birational morphism, $\gamma\circ h(R')$ is a line
     (plus possibly some isolated points).
     This way there exists a smooth rational curve $\beta\subset B'$
     such that $$K_V\widetilde{\beta}=1,$$ where $\widetilde{\beta}\subset
     R'$ is the support of $\pi^*(\beta).$ The adjunction formula
     gives $\widetilde{\beta}^2=-3,$ thus $\beta^2=-6.$
     Notice that $\widetilde{\beta}$ is the only component of $R'$ which is not
     contracted by the map $\gamma\circ h$.

     Let $$\alpha:=B'-\beta\subset W,$$
     $$\overline{\beta}:=\rho(\beta),\ \overline{\alpha}:=\rho(\alpha)\subset P.$$
     When $\alpha$ is non-empty, the support of $\pi^*(\alpha)$
     is an union of $(-2)$-curves,
     since it is contracted by $\gamma\circ h.$ Equivalently $\alpha$
     is a disjoint union of $(-4)$-curves.
     \\\\
     {\bf Claim 5:} {\em We have $$K_W^2\geq -2.$$}

     \noindent {\em Proof}\ : Consider the Chern number $c_2$ and the
     second Betti number $b_2.$ It is well known that, for a surface $X,$
     $$c_2(X)=12\chi(\mathcal{O}_X)-K_X^2,\ \ b_2(X)=c_2(X)-2+4q(X).$$
     Therefore $$b_2(W)=22-K_W^2,\ \ b_2(V)=b_2(S)+t=11+13=24.$$
     The inequality $K_W^2\geq -2$ follows from the fact $b_2(V)\geq b_2(W).$ $\diamond$
     \\

     From Claim 5, we conclude that the resolution of $\overline{B'}$ blows-up
     at most two double points, thus $$B'^2\geq -2+2(-4)=-10=\beta^2+(-4).$$
     This implies that $\alpha$ is a smooth $(-4)$-curve when $\alpha\ne 0$.
     \\\\
     {\bf Claim 6:} {\em Only the following possibilities can occur:
     \begin{itemize}
       \item [$\cdot$] $\overline{\beta}$ has one double point and no other
       singularity, or
       \item [$\cdot$] $\overline{\alpha},$ $\overline{\beta}$ are
       smooth, \ $\overline{\alpha}\overline{\beta}=2.$
     \end{itemize}}

     \noindent {\em Proof}\ : Recall that $B'=\alpha+\beta$ is
     contained in fibres of $f_A$ and, since $\overline{B'}$ has
     only negligible singularities, then also
     $\overline{B'}=\overline{\alpha}+\overline{\beta}$ is contained
     in fibres of $\overline{f_A}.$ In particular
     $\overline{\alpha}^2,\overline{\beta}^2\leq 0.$

     If $\overline{\alpha}$ is singular, then it has arithmetic genus
     $p_a(\overline{\alpha})=1$ and $\overline{\alpha}^2=0.$ But
     then $\overline{\alpha}$ has the same support of a fibre of
     $\overline{f_A},$ which is a contradiction because $\overline{f_A}$
     is not elliptic. Therefore $\overline{\alpha}$ is smooth.

     Since $K_P\overline{\alpha}\geq 0,$ $K_P\overline{B'}=2$ implies $K_P\overline{\beta}\leq 2.$
     We know that $\beta$ is a smooth rational curve and
     $\beta^2=-6,$ thus $K_W\beta=4.$ If $\overline{\beta}$ is smooth, then
     one must have $\overline{\alpha}\overline{\beta}>1.$ From Claim 5 the only
     possibility in this case is $\overline{\alpha}\overline{\beta}=2.$
     If $\overline{\beta}$ is singular, then
     $\overline{\beta}^2\leq 0$ implies that
     $\overline{\beta}$ has one ordinary double point and no other
     singularity. $\diamond$
     \\

     Let $D:=\overline{\beta}$ if $\overline{\beta}$ is singular.
     Otherwise let $D:=\overline{\alpha}+\overline{\beta}.$

     The 2-connected divisor
     $\widetilde{D}:=\frac{1}{2}(\rho\circ\pi)^*(D)$ has arithmetic genus
     $p_a(\widetilde{D})=1.$
     We know that $(K_V+n\Lambda)\widetilde{D}=1$ (because $K_VR'=1$) and that
     $\widetilde{D}$ contains a component $A$ such that $(K_V+n\Lambda)A=0$
     (because $D$ has at least one negligible singularity).
     These two facts imply, from \cite[Proposition A.5, (ii)]{CFM},
     that the relative canonical map $\gamma$ has a base point in $\widetilde{D}.$
     As mentioned above, $\gamma$ is a morphism, which is a contradiction.\\

Finally the assertion about ${\rm deg}(\phi_2)$ in Case 2): we have
proved that $S$ has a genus 2 fibration, so it has an
hyperelliptic involution $j.$ The bicanonical map $\phi_2$ factors
through both $i$ and $j,$ thus ${\rm deg}(\phi_2)\geq 4.$

This finishes the proof of case a),\ (ii) of Theorem \ref{next}.
\\\\
We end the proof of Theorem \ref{next} with {\bf Case a),\ (iii)}: A
surface of general type with a genus 2 fibration and $p_g=q=1$
satisfies $K^2\leq 6$ (see \cite{X1}). Denote by $j$ the map such
that $\phi_2=j\circ i.$ The quotient $S/i$ is a $K3$ surface thus,
from \cite{S}, ${\rm deg}(j)\leq 2.$ Analogously to Case 2, ${\rm
deg}\ (\phi_2)\geq 4,$ thus ${\rm deg}(j)=2,$ ${\rm deg}(\phi_2)=4$
and then $K_S^2\ne 2$ (see Case 1).

It follows from \cite[p. 66]{X1} that, if the genus 2 fibration of
$S$ has a rational basis, then $K_S^2=3.$ It is shown in \cite{P}
that, in these conditions, ${\rm deg}(\phi_2)=2.$ We then conclude
that the genus 2 fibration of $S$ is the Albanese fibration.
\\\\
Examples for case a),\ (iii) with $K_S^2=3,4$ or 5 were given by
Catanese in \cite{C2}. The existence of the other cases is proved in
the next section.\Qed
\section{Examples}\label{Examples}
In this section we will construct smooth minimal surfaces of general
type $S$ with $p_g(S)=q(S)=1$ having an involution $i$ such that the bicanonical
map $\phi_2$ of $S$ is composed with $i$ and:
\begin{description}
   \item[1)] $K_S^2=6,$ $g=3$, ${\rm deg}(\phi_2)=2,$ $S/i$ is birational to a $K3$ surface;
   \item[2)] $K_S^2=4,$ $g=2,$ ${\rm Kod}(S/i)=1;$
   \item[3)] $K_S^2=2,$ $g=2,$ ${\rm Kod}(S/i)=2,$
\end{description}
where $g$ denotes the genus of the Albanese fibration of $S.$\\\\
{\bf Example 1:}\\\\ In \cite{T} Todorov gives the following
construction of a surface of general type $S$ with $p_g(S)=1,$
$q(S)=0$ and $K_S^2=8.$ Consider a Kummer surface $Q$ in $\P^3,$
i.~e. a quartic having has only singularities 16 nodes (ordinary
double points). Let $G\subset Q$ be the intersection of $Q$ with a
general quadric, $\widetilde{Q}$ be the minimal resolution of $Q$
and $\widetilde{G}\subset \widetilde{Q}$ be the pullback of $G.$ The
surface $S$ is the minimal model of the double cover
$\pi:V\rightarrow\widetilde{Q}$ ramified over
$\widetilde{G}+\sum_1^{16} A_i,$ where $A_i\subset \widetilde{Q},$
$i=1,\ldots,16,$ are the $(-2)$-curves which contract to the nodes
of $Q.$

It follows from the double cover formulas (cf. \cite[Ch. V, \S 22]{BPV})
that the imposition of a quadruple point to the branch locus
decreases $K^2$ by 2 and the Euler characteristic $\chi$ by 1.

We will see that we can impose a quadruple point to the branch locus
of the Todorov construction, thus obtaining $S$ with $K_S^2=6.$ In
this case I claim that $p_g(S)=q(S)=1.$ In fact, let $W$ be the surface $\widetilde Q$
blown-up at the quadruple point, $E$ be the corresponding $(-1)$-curve,
$B$ be the branch locus and $L$ be the line bundle such that $2L\equiv B.$
From formula (\ref{Mrg2}) in Section  \ref{Involution}, one has 
$h^0(W,\mathcal O_W(2E+L))=0$ (thus the bicanonical map of $V$ factors through $\pi$), 
hence also $h^0(W,\mathcal O_W(E+L))=0$
and then $$p_g(S)=p_g(W)+h^0(W,\mathcal O_W(E+L))=1.$$

We will see that ${\rm deg}(\phi_2)=2,$
hence $\phi_2(S)$ is a $K3$ surface and so $S$ has no genus 2
fibration.
\\

First we need to obtain an equation of a Kummer surface. The
Computational Algebra System {\em Magma} has a direct way to do
this, but I prefer to do it using a beautiful construction that I learned
from Miles Reid.

We want a quartic surface $Q\in\P^3$ whose singularities are exactly
16 nodes. Projecting from one of the nodes to $\P^2,$ one realizes
the "Kummer" surface as a double cover
$$\psi:X\longrightarrow\P^2$$
with branch locus the union of 6 lines $L_i$ (see \cite[p.
774]{GH}), each one tangent to a conic $C$ (the image of the
projection point) at a point $p_i.$ The surface $X$ contains 15
nodes (from the intersection of the lines) and two $(-2)$-curves
(the pullback $\psi^*(C)$) disjoint from these nodes. To obtain a
Kummer surface we have just to contract one of these curves.

Denote also by $L_i$ the defining polynomial of each line $L_i.$ An
equation for $X$ is $z^2=L_1\cdots L_6$ in the weighted projective
space $\P(3,1,1,1),$ with coordinates $(z,x_1,x_2,x_3)$. We will see
that this equation can be written in the form $AB+DE=0,$ where the
system $A=B=D=E=0$ has only the trivial solution and $B,E$ are the
defining polynomials of one of the $(-2)$-curves in $\psi^*(C).$ Now
consider the surface $X'$ given by $Bs=D,\ Es=-A$ in the space
$\P(3,1,1,1,1)$ with coordinates $(z,s,x_1,x_2,x_3).$ There is a
morphism $X\rightarrow X'$ which restricts to an isomorphism
$$X-\{B=E=0\}\longrightarrow X'-\{[0:1:0:0:0]\}$$ and which
contracts the curve $\{B=E=0\}$ to the point $[0:1:0:0:0].$ This is
an example of {\em unprojection} (see \cite{R}).

The variable $z$ appears isolated in the equations of $X',$
therefore eliminating $z$  we obtain the equation of the Kummer $Q$
in $\P^3$ with variables $(s,x_1,x_2,x_3).$ All this calculations
will be done using Magma.

In what follows a line preceded by $>$ is an input line, something
preceded by // is a comment. A $\backslash$ at the end of a line
means continuation in the next line. The other lines are output
ones.

\

\begin{verbatim}
> K<e>:=CyclotomicField(6);//e denotes the 6th root of unity.
> //We choose a conic C with equation x1x3-x2^2=0 and fix the
> //p_i's: (1:1:1), (e^2:e:1), (e^4:e^2:1), (e^6:e^3:1),
> //(e^8:e^4:1), (e^10:e^5:1).
> R<z,s,x1,x2,x3>:=PolynomialRing(K,[3,1,1,1,1]);
> g:=&*[e^(2*i)*x1-2*e^i*x2+x3:i in [0..5]];
> //g is the product of the defining polynomials
> //of the tangent lines L_i to C at p_i.
> X:=z^2-g;
> X eq (z+x1^3-x3^3)*(z-x1^3+x3^3)+4*(x1*x3-4*x2^2)^2*\
> (-x1*x3+x2^2);//The decomposition AB+DE.
true
> i:=Ideal([s*(z-x1^3+x3^3)-4*(x1*x3-4*x2^2)^2,\
> s*(x1*x3-x2^2)-(z+x1^3-x3^3)]);
> j:=EliminationIdeal(i,1);
> j;
Ideal of Graded Polynomial ring of rank 5 over K
Lexicographical Order Variables: z, s, x1, x2, x3
Variable weights: 3 1 1 1 1 Basis:
[-1/2*s^2*x1*x3+1/2*s^2*x2^2+s*x1^3-s*x3^3+2*x1^2*x3^2-
16*x1*x2^2*x3+32*x2^4]
> 2*Basis(j)[1];
-s^2*x1*x3+s^2*x2^2+2*s*x1^3-2*s*x3^3+4*x1^2*x3^2-
32*x1*x2^2*x3+64*x2^4
> //This is the equation of the Kummer Q.
\end{verbatim}

\

We want to find a quadric $H$ such that $H\bigcap Q$ is a
reduced curve $\overline{B'}$ having an ordinary quadruple point $pt$ as only
singularity. Since the computer is not fast enough while working with more
than 5 or 6 variables, we first need to think what the most
probable case is.

Like we have seen in Section \ref{q1}, the branch locus
$B'+\sum_1^{16}A_i$ is contained in 4 fibres $F_A^1,\ldots,F_A^4$ of
a fibration $f_A$ of $W,$ where $W$ is the resolution of $Q$
blown-up at $pt$ and the $A_i$'s are the $(-2)$-curves which
contract to the nodes of $Q.$

Of course we have a quadric intersecting $Q$ at a curve with a
quadruple point $pt:$ the tangent space $T$ to $Q$ at $pt$ counted
twice. But this one is double, so we need to find an irreducible one
(and these two induce $f_A$), the curve $\overline{B'}.$ These
curves $2T$ and $\overline{B'}$ are good candidates for
$\overline{F_A^1}$ and $\overline{F_A^2}$ (in the notation of
Sections \ref{q1} and \ref{Classification}). If this configuration
exists, then the 16 nodes must be contained in the other two fibres,
$\overline{F_A^3}$ and $\overline{F_A^4}.$ These fibres are
divisible by 2, because $\overline{F_A^1}=2T,$ and are double
outside the nodes.
Since in a $K3$ surface only 0, 8 or 16 nodes can have sum divisible by 2,
it is reasonable to try the following configuration: each of $\overline{F_A^3}$ and
$\overline{F_A^4}$ contain 8 nodes with sum divisible by 2 and is double outside the nodes.

It is well known (see e.~g. \cite{GH}) that the Kummer surface $Q$
has 16 double hyperplane sections $T_i$ such that each one contains
6 nodes of $Q$ and that any two of them intersect in 2 nodes.
The sum of the 8 nodes contained in
$$N:=(T_1\cup T_2)\backslash(T1\cap T_2)$$ is divisible by 2. Magma
will give 3 generators $h_1,$ $h_2,$ $h_3$ for the linear system of
quadrics through these nodes.

\

\begin{verbatim}
> K<e>:=CyclotomicField(6);
> P3<s,x1,x2,x3>:=ProjectiveSpace(K,3);
> F:=-s^2*x1*x3+s^2*x2^2+2*s*x1^3-2*s*x3^3+4*x1^2*x3^2-\
> 32*x1*x2^2*x3+64*x2^4;
> Q:=Scheme(P3,F);/*The Kummer*/; SQ:=SingularSubscheme(Q);
> T1:=Scheme(P3,x1-2*x2+x3); T2:=Scheme(P3,s);
> N:=Difference((T1 join T2) meet SQ, T1 meet T2);
> s:=SetToSequence(RationalPoints(N));
> //s is the sequence of the 8 nodes.
> L:=LinearSystem(P3,2);
>//This will give the h_i's:
> LinearSystem(L,[P3!s[i] : i in [1..8]]);
Linear system on Projective Space of dimension 3
  Variables: s, x1, x2, x3 with 3 sections:

s*x1-2*x1^2-4*x1*x2-2*x1*x3+8*x2^2+4*x2*x3+2*x3^2
s*x2-2*x1^2-4*x1*x2+4*x2*x3+2*x3^2
s*x3-2*x1^2-4*x1*x2+2*x1*x3-8*x2^2+4*x2*x3+2*x3^2
\end{verbatim}

\

Now we want to find a quadric $H$ in the form $h_1+bh_2+ch_3,$
for some $b,c$ (or, less probably, in the form $bh_2+ch_3$) such
that the projection of $H\bigcap Q$ to $\P^2$ (by elimination) is
a curve with a quadruple point. To find a quadruple point we just
have to impose the annulation of the derivatives up to order 3 and
ask Magma to do the rest.

\

\begin{verbatim}
> R<s,b,c,x1,x2,x3>:=PolynomialRing(Rationals(),6);
> F:=-s^2*x1*x3+s^2*x2^2+2*s*x1^3-2*s*x3^3+4*x1^2*x3^2-\
> 32*x1*x2^2*x3+64*x2^4;
> h1:=s*x1-2*x1^2-4*x1*x2-2*x1*x3+8*x2^2+4*x2*x3+2*x3^2;
> h2:=s*x2-2*x1^2-4*x1*x2+4*x2*x3+2*x3^2;
> h3:=s*x3-2*x1^2-4*x1*x2+2*x1*x3-8*x2^2+4*x2*x3+2*x3^2;
> H:=h1+b*h2+c*h3;
> I:=ideal<R|[F,H]>;
> I1:=EliminationIdeal(I,1);
> q0:=Evaluate(Basis(I1)[1],x3,1);//We work in the affine plane.
> R4<B,C,X1,X2>:=PolynomialRing(Rationals(),4);
> h:=hom<R->R4|[0,B,C,X1,X2,0]>;
> q:=h(q0);q1:=Derivative(q,X1);q2:=Derivative(q,X2);
> q3:=Derivative(q1,X1);q4:=Derivative(q1,X2);q5:=Derivative\
> (q2,X2);q6:=Derivative(q3,X1);q7:=Derivative(q3,X2);
> q8:=Derivative(q4,X2);q9:=Derivative(q5,X2);
> A4:=AffineSpace(R4);
> S:=Scheme(A4,[q,q1,q2,q3,q4,q5,q6,q7,q8,q9]);
> Dimension(S);
0
> PointsOverSplittingField(S);
\end{verbatim}

\

This last command gives the points of $S,$ as well as the necessary
equations to define the field extensions where they belong.
There are various solutions. One of them gives the desired quadruple point.
The confirmation is as follows:

\

\begin{verbatim}
> R<x>:=PolynomialRing(Rationals());
> K<r13>:=ext<Rationals()|x^4 + x^3 + 1/4*x^2 + 3/32>;
> P3<s,x1,x2,x3>:=ProjectiveSpace(K,3);
> F:=-s^2*x1*x3+s^2*x2^2+2*s*x1^3-2*s*x3^3+4*x1^2*x3^2-\
> 32*x1*x2^2*x3+64*x2^4;
> b:=64/55*r13^3-272/55*r13^2-96/55*r13-46/55;
> c:=-2176/605*r13^3+448/605*r13^2+624/605*r13-361/605;
> H:=(s*x1-2*x1^2-4*x1*x2-2*x1*x3+8*x2^2+4*x2*x3+2*x3^2)+\
> b*(s*x2-2*x1^2-4*x1*x2+4*x2*x3+2*x3^2)+\
> c*(s*x3-2*x1^2-4*x1*x2+2*x1*x3-8*x2^2+4*x2*x3+2*x3^2);
> Q:=Scheme(P3,F);
> C:=Scheme(Q,H);
> IsReduced(C);
false
> RC:=ReducedSubscheme(C);
> #SingularPoints(RC);//# means "number of".
1
> HasSingularPointsOverExtension(RC);
false
> pt:=Representative(SingularPoints(RC));
> pt in SingularSubscheme(Q);//pt is not a node of Q.
false
> T:=DefiningPolynomial(TangentSpace(Q,pt));
> T2:=Scheme(Q,T^2);
> #RationalPoints(T2 meet C);
1
> pt in RationalPoints(T2 meet C);
true
> HasPointsOverExtension(T2 meet C);
false
\end{verbatim}
This way $T2$ and $C$ generate a pencil with a quadruple base point
and the curve $\overline{B'}$ is a general element of this pencil.

\

Finally, it remains to be shown that the degree of the bicanonical map
$\phi_2$ is 2. As $(2K_S)^2=24,$ it suffices to show that
$\phi_2(S)$ is of degree 12. Since, in the notation of diagram
(\ref{CmtDgr}), $h^*|2K_S|=\pi^*|2K_W+B'|$ then $\phi_2(S)$ is the
image of $W$ via the map $\tau:W\rightarrow \phi_2(S)$ given by
$|2K_W+B'|.$ The projection of this linear system on $Q$ is
the linear system of the quadrics whose intersection with $Q$ has a
double point at $pt.$ In order to easily write this linear system,
we will translate the point $pt$ to the origin (in affine
coordinates).

\

\begin{verbatim}
> QA:=AffinePatch(Scheme(P3,F),4);
> p:=Representative(RationalPoints(AffinePatch(Cluster(pt),4)));
> A3<x,y,z>:=Ambient(QA);
> psi:=map<A3->A3|[x-p[1],y-p[2],z-p[3]]>;Q0:=psi(QA);
> FA:=DefiningPolynomial(Q0);
> j:=[Evaluate(Derivative(FA,A3.i),Origin(A3)):i in [1,2,3]];
> J:=LinearSystem(A3,[j[1]*x+j[2]*y+j[3]*z,x^2,x*y,x*z,y^2,y*z,\
> z^2]);
> P6:=ProjectiveSpace(K,6);
> tau:=map<A3->P6|Sections(J)>;
> Degree(tau(Q0));
12
\end{verbatim}
{\bf Example 2:}
\\\\
Here we will construct a surface of general type $S,$ with $p_g=q=1$
and $K^2=4,$ as the minimal model of a double cover of a surface $W$
such that ${\rm Kod} (W)=p_g(W)=1$ and $q(W)=0.$
\\\\
{\bf Step 1:} {\em Construction of $W.$}
\\
Consider five distinct lines $L_1,\ldots,L_5\subset\P^2$ meeting in
one point $p_0.$ Let $p_1\in L_4,$ $p_2,p_3\in L_5$ be points
distinct from $p_0.$ Choose three distinct non-degenerate conics,
$C_1,C_2,C_3,$ tangent to $L_4$ at $p_1$ and passing through
$p_2,p_3.$
Define $$D:=L_1+\ldots+L_4+C_1+C_2.$$ Denote by $p_4,\ldots,p_{15}$
the 12 nodes of $D$ contained in $L_1+L_2+L_3.$ To resolve the
$(3,3)$ point of $D$ at $p_1$ we must do two blow-ups: one at $p_1$
and other at an infinitely near point $p_1'.$ Let
$\mu:X\rightarrow\P^2$ be the blow-up with centers
$p_0,p_1,p_1',p_2,\ldots,p_{15}$ and
$E_0,E_1,E_1',E_2,\ldots,E_{15}$ be the corresponding exceptional
divisors (with self-intersection $-1$).
Consider
$$D':=\mu^*(D)-4E_0-2E_1-4E_1'-2\sum_2^{15}E_i.$$
Let $\psi:\widetilde{X}\rightarrow X$ be the double cover of $X$
with branch locus $D'.$ The surface $\widetilde{X}$ is the canonical
resolution of the double cover of $\P^2$ ramified over $D.$ Let $W$
be the minimal model of $\widetilde{X}$ and $\nu$ be the
corresponding morphism.
$$
\begin{CD}\ \widetilde{X}@>\nu>>W\\ @V\psi VV \\ X@>\mu >> \P^2
\end{CD}
$$
Notice that $\nu$ contracts two $(-1)$-curves contained in
$(\mu\circ\psi)^*(L_4).$

We have $K_X\equiv -\mu^*(3L)+E_1'+\sum_0^{15}E_i,$ where $L$
denotes a general line of $\P^2.$ Hence, using the double cover
formulas (cf. (\ref{relate})),
$$K_{\widetilde{X}}\equiv\psi^*\left(K_X+\frac{1}{2}D'\right)
\equiv\psi^*(\mu^*(L)-E_0-E_1')\equiv\psi^*(\widehat{L_4}+(E_1-E_1')+E_1'),$$
where $\widehat{L_4}\subset X$ is the strict transform of $L_4.$
Since $\widehat{L_4}$ and $E_1-E_1'$ are $(-2)$-curves contained in
the branch locus $D',$ then $\frac{1}{2}\psi^*(\widehat{L_4})$ and
$\frac{1}{2}\psi^*(E_1-E_1')$ are $(-1)$-curves in $\widetilde{X},$
thus
$$K_W\equiv \nu(\psi^*(E_1')).$$
The divisor $2\nu(\psi^*(E_1'))\equiv 2K_W$ is a (double) fibre of
the elliptic fibration of $W$ induced by the pencil of lines through
$p_0.$ So $p_g(W)=1$ and $W$ has Kodaira dimension 1.

From (\ref{relate}) one has
$$\chi(\mathcal{O}_W)=2+\frac{1}{8}D'(2K_X+D')=2+\frac{1}{8}(28-28)=2.$$
\\\\
{\bf Step 2:} {\em The branch locus in $W.$}
\\
Since the strict transforms
$\widehat{L_1},\ldots,\widehat{L_4}\subset X$ are in the branch
locus $D',$ then there are curves $l_1,\ldots,l_4\subset
\widetilde{X}$ such that
$$(\mu\circ\psi)^*(L_1+\cdots+L_4)=2l_1+\cdots+
2l_4+4\psi^*(E_0)+\psi^*(E_1-E_1')+2\psi^*(E_1')+\sum_4^{15}A_i,$$
where each $A_i:=\psi^*(E_i)$ is a $(-2)$-curve. But also $E_1-E_1'$
is in the branch locus, thus $\psi^*(E_1-E_1')\equiv 0\ ({\rm mod\
2})$ and then
$$\sum_4^{15}A_i\equiv 0\ ({\rm mod\ 2}).$$
The strict transform $\widehat{L_5}$ is a $(-2)$-curve which do not
intersect $D'$ thus
$$\psi^*(\widehat{L_5})=A_{16}+A_{17},$$ with $A_{16},A_{17}$
disjoint $(-2)$-curves.

Denote by $\widehat{C_3}\subset X$ the strict transform of the conic
$C_3.$ We have
$$(\mu\circ\psi)^*(C_3+L_4+L_5)=$$$$=\psi^*(\widehat{C_3})+2l_4+A_{16}+A_{17}
+2\psi^*(E_0+\cdots+E_3)+2\psi^*(E_1')\equiv 0\ (\rm mod\ 2).$$
With this we conclude that $$\psi^*(\widehat{C_3})+\sum_4^{17}
A_i\equiv 0\ ({\rm mod\ 2}).$$ Notice that
$F\cdot\nu\left(\psi^*(\widehat{C_3})\right)=4$ for a fibre $F$ of
the elliptic fibration of $W,$ thus
$K_W\cdot\nu\left(\psi^*(\widehat{C_3})\right)=2.$
\\\\
{\bf Step 3:} {\em Construction of $S.$}
\\
Let $\pi:V\rightarrow W$ be the double cover with branch locus
$$B:=\nu\left(\psi^*(\widehat{C_3})+\sum_4^{17} A_i\right)$$
and $S$ be the minimal model of $V.$ From the double cover formulas
(\ref{relate}) we obtain
$$2K_V^2=(2K_W+B)^2=4K_W^2+4K_WB+B^2=4\cdot 0+4\cdot 2+(-28)=-20$$
and, by contraction of the $(-1)$-curves
$\frac{1}{2}\pi^*(\nu(A_i)),$
$$K_S^2=K_V^2+14=-10+14=4.$$

Let $L:\equiv\frac{1}{2}B.$ Formulas (\ref{relate}) give
$$\chi(\mathcal{O}_S)=2\chi(\mathcal{O}_W)+\frac{1}{2}L(K_W+L)=4-3=1.$$
Using now formula (\ref{Mrg2}) we obtain
$h^0(W,\mathcal{O}_W(2K_W+L))=0,$ which means that the bicanonical
map of $V$ factors through $\pi.$

Because $K_W$ is effective then also
$h^0(W,\mathcal{O}_W(K_W+L))=0$ and
$$p_g(S)=p_g(W)+h^0(W,\mathcal{O}_W(K_W+L))=1.$$

Hence $q(S)=1$ and then, as we noticed in the beginning of Section \ref{Classification},
the curve $\nu\left(\psi^*(\widehat{C_3})\right)$ is contained in the
fibration of $W$ which induces the Albanese fibration of $S.$ As
$\nu\left(\psi^*(\widehat{C_3})\right)^2=0,$ we conclude that the
Albanese fibration of $S$ is the one induced by the pencil
$|\widehat{C_3}|.$ It is of genus 2 because
$\widehat{C_3}D'=\widehat{C_3}(\widehat{L_1}+\widehat{L_2}+\widehat{L_3})=6.$
\\\\
{\bf Example 3:}
\\\\
Now we will obtain a surface of general type $S,$ with $p_g=q=1$ and
$K^2=2,$ as the minimal model of a double cover of a surface of
general type $W$ such that $K_W^2=p_g(W)=1$ and $q(W)=0.$
\\\\
{\bf Step 1:} {\em Construction of $W.$}\\
Let $p_0,\ldots,p_3\in\P^2$ be distinct points and $L_i$ be the line
through $p_0$ and $p_i,$ $i=1,2,3.$ For each $j\in\{1,2,3\}$ let
$C_j$ be the conic through $p_1,p_2,p_3$ tangent to the $L_i$'s
except for $L_j.$ Denote by $D$ a general element of the linear
system generated by $3C_1+2L_1,$ $3C_2+2L_2$ and $3C_3+2L_3.$ The
singularities of $D$ are a $(3,3)$-point at $p_i,$ tangent to $L_i,$
$i=1,2,3,$ and a double point at $p_0.$ Let $L_4$ be a line through
$p_0$ transverse to $D.$

Denote by $W'$ the canonical resolution of the double cover of
$\P^2$ with branch locus $$D+L_1+\ldots+L_4$$ and by $W$ the minimal
model of $W'.$ The formulas of \cite[Ch. V, \S 22]{BPV} give
$\chi(W)=2$ and $K_W^2=1$ (notice that the map $W'\rightarrow W$
contracts three $(-1)$-curves contained in the pullback of
$L_1+L_2+L_3$). Since $K^2\geq 2p_g$ for an irregular surface
(\cite{De}), $W$ is regular and then $p_g(W)=\chi(W)-1=1.$
\\\\
{\bf Step 2:} {\em The branch locus in $W.$}\\
The pencil of lines through $p_0$ induces a (genus 2) fibration of $W.$ Let
$F_i$ be the fibre induced by $L_i,$ $i=1,\ldots,4.$ The fibre $F_4$
is the union of six disjoint $(-2)$-curves (corresponding to the
nodes of $D-p_0$) with a double component (the strict transform of
$L_4$). Each $F_i,$ $i=1,2,3,$ is the union of two $(-2)$-curves
with a double component (cf. \cite[\S 1]{X1}). Thus $F_1+\cdots+F_4$
contain disjoint $(-2)$-curves $A_1,\ldots,A_{12}$ such that
$$\sum_1^{12}A_i\equiv 0\ ({\rm mod}\ 2).$$
\\\\
{\bf Step 3:} {\em Construction of $S.$}\\
Let $V$ be the double cover of $W$ with branch locus
$\sum_1^{12}A_i$ and $S$ be the minimal model of $V.$ From
(\ref{relate}) we obtain $\chi(\mathcal{O}_S)=1$ and $K_V^2=-10.$
The $A_i$'s lift to $(-1)$-curves in $V,$ thus $K_S^2=-10+12=2.$ We
have $1=p_g(W)\leq p_g(S),$ hence $q(S)\ne 0$ and then $2=K_S^2\geq
2p_g(S).$ So $p_g(S)=q(S)=1.$

The genus 2 fibration of $W$ induces the Albanese fibration of $S.$

\bigskip

\bigskip

\medskip

\noindent Carlos Rito
\\ Departamento de Matem\' atica
\\ Universidade de Tr\' as-os-Montes e Alto Douro
\\ 5000-911 Vila Real
\\ Portugal
\\\\
\noindent {\it e-mail}: crito@utad.pt
\end{document}